\documentclass[12pt]{amsart}
\usepackage{amsmath, amssymb,graphicx,caption, subcaption, mathtools, enumitem}
\usepackage[pdftex,hyperindex=true]{hyperref}
\usepackage{url}
\usepackage{etoolbox}
\patchcmd{\section}{\scshape}{\bfseries}{}{}
\patchcmd{\subsection}{\bfseries}{\itshape}{}{}

\makeatletter
\def\@seccntformat#1{%
  \protect\textup{\protect\@secnumfont
    \ifnum\pdfstrcmp{section}{#1}=0 \bfseries\fi
    \ifnum\pdfstrcmp{subsection}{#1}=0 \itshape \fi
    \csname the#1\endcsname
    \protect\@secnumpunct
  }%
}  
\makeatother

\usepackage[utf8]{inputenc}
\usepackage[english]{babel}

\theoremstyle{plain}

\newtheorem{theorem}{Theorem}[section]

\newtheorem{conjecture}{Conjecture}[section]

\theoremstyle{definition}

\newtheoremstyle{note}
{3pt}
{3pt}
{\itshape}
{}
{\itshape}
{:}
{.5em}
{}

\theoremstyle{note}
\newtheorem{remark}[theorem]{Remark}

\theoremstyle{plain} 
\newcommand{\thistheoremname}{}
\newtheorem*{genericthm}{\thistheoremname}

\numberwithin{equation}{section}

\newcommand{\acknowledge}{\subsection*{Acknowledgements}}

\def\Dbar{\leavevmode\lower.6ex\hbox to 0pt{\hskip-.23ex \accent"16\hss}D}

\begin{document}

\title{Towards the Atiyah-Sutcliffe conjectures for coplanar hyperbolic points}
\author{J. Malkoun}
\address{Department of Mathematics and Statistics\\
Faculty of Natural and Applied Sciences\\
Notre Dame University-Louaize, Zouk Mosbeh\\
P.O.Box: 72, Zouk Mikael,
Lebanon}
\email{joseph.malkoun@ndu.edu.lb}

\date{September 2, 2019}

\maketitle

\begin{abstract} The Atiyah-Sutcliffe normalized determinant function $D$ is a smooth complex-valued function on $C_n(H^3)$, where 
$C_n(H^3)$ denotes the configuration space of $n$ distinct points in hyperbolic $3$-space $H^3$. The hyperbolic version of the Atiyah-Sutcliffe conjecture 
$1$ (AS conjecture $1$) states that $D$ is nowhere vanishing. AS conjecture $2$ (hyperbolic version) is the stronger statement that $|D(\mathbf{x})| \geq 1$ 
for any $\mathbf{x} \in C_n(H^3)$.

In this short article, we prove AS conjecture $2$ for hyperbolic convex coplanar quadrilaterals, that is for configurations of $4$ points in $H^2$ with none of 
the points in the configuration lying in the convex hull of the other three. We also obtain Y. Zhang and J. Ma's result, namely AS conjecture $1$ for non-convex 
quadrilaterals in $H^2$. Finally, we find an explicit lower bound for $|D|$ depending on $n$ only for the natural ``star-based'' variant of the AS problem, for convex 
coplanar hyperbolic configurations. The latter result holds for any $n \geq 2$.

The proofs for $n=4$ make use of the symbolic library of Python. The proof of the general result follows from a general formula for the determinant. In all 
these cases, $D$ can be expanded as a linear combination of non-negative rational functions with positive coefficients. 
\end{abstract}

\section{Introduction} \label{intro}

Denote by $C_n(H^3)$ the configuration space of $n$ distinct points $\mathbf{x}_1,\ldots$, $\mathbf{x}_n$ in $H^3$. For short, we denote 
such a configuration by
\[ \mathbf{x} = (\mathbf{x}_1,\ldots,\mathbf{x}_n) \in C_n(H^3) .\]
We use the open unit $3$-ball (centered at the origin) Poincar\'{e} model of $H^3$. The boundary sphere $\partial H^3$ will be referred to as 
the sphere at infinity, and denoted by $S^2_{\infty}$.

Given a configuration $\mathbf{x} \in C_n(H^3)$, we shall define $n$ complex polynomials $p_a(t)$, $1 \leq a \leq n$, of degree at most $n-1$ in 
one complex variable $t$, with each polynomial only defined up to scaling by a non-zero complex factor. 

Before defining the polynomials $p_a(t)$, we define the stereographic projection $s: S^2_{\infty} \to \hat{\mathbb{C}}$ by the expression
\[ s(x,y,z) = \left\{ \begin{array}{ll} \frac{x+iy}{1-z}, &\text{ if $(x,y,z) \neq (0,0,1)$} \\
                          \infty, &\text{ if $(x,y,z) = (0,0,1)$.} \end{array} \right. \]
Stereographic projection is a diffeomorphism, and allows us transport the complex structure from the Riemann sphere 
$\hat{\mathbb{C}} = P^1(\mathbb{C})$ to the $2$-sphere $S^2_{\infty}$, and thus describe the location of each point on the sphere 
at infinity using a \emph{complex} coordinate, which is allowed to take the value $\infty$.

\emph{From now on, we shall identify $S^2_{\infty}$ with the Riemann sphere $\hat{\mathbb{C}}$ via $s$.}

Given a pair of indices $(a,b)$, with $1 \leq a,b \leq n$ and $a \neq b$, we define the point $t_{ab} \in S^2_{\infty}$ to be the limiting 
point of the hyperbolic ray starting at $\mathbf{x}_a$ and passing through $\mathbf{x}_b$. We will refer to the $t_{ab}$ as ``directions'', in analogy 
with the Euclidean version of the problem (see for instance \cite{Atiyah2000}, \cite{At1} or \cite{AS2}). Thus, by our previous remark, each 
$t_{ab}$ is a complex number, or possibly infinity.

We are now ready to define $p_a(t)$, for $1 \leq a \leq n$, as the complex polynomial of degree at most $n-1$, having as roots the $t_{ab}$, 
for $1 \leq b \leq n$ and $b \neq a$. So
\begin{equation} p_a(t) = \prod_{b\neq a} (t - t_{ab}) . \label{observer-based} \end{equation}
This expression makes sense if none of the $t_{ab}$ is infinite. How do we interpret a factor such as $t - t_{ab}$ if $t_{ab} = \infty$? The key to 
resolving this issue is to use homogeneous coordinates on $P^1(\mathbb{C})$. The difference gets then replaced by a $2$-by-$2$ complex 
determinant. We thus see that if $t_{ab} = \infty$, then the factor $t - t_{ab}$ can be replaced by a non-zero constant, such as $1$, say.

\begin{conjecture}[AS conjecture $1$] (AS stands for Atiyah-Sutcliffe) states that no matter which $\mathbf{x} \in C_n(H^3)$ one starts with, the associated $n$ 
polynomials $p_a(t)$, $1 \leq a \leq n$, are linearly independent over $\mathbb{C}$. \label{AS1} \end{conjecture}

It will be convenient to introduce the following notation
\[ p_{ab}(t) = t - t_{ab} .\]

The normalized AS determinant function $D: C_n(H^3) \to \mathbb{C}$ can now be defined.
\[ D(\mathbf{x}) = \frac{\det(p_1,\ldots,p_n) }{\prod_{a<b} \det(p_{ab},p_{ba})} \]

In the previous expression, each $p_a$ is identified with its $n$-dimensional coefficient vector, corresponding to increasing powers of $t$, and 
the numerator is thus the determinant of the $n$-by-$n$ complex matrix having $p_b$ as its $b$-th column. Similarly, each factor in the denominator 
is a complex $2$-by-$2$ determinant, constructed similarly (each $p_{ab}$ has degree at most $1$ in $t$, so that its vector of coefficients is $2$-dimensional).

\begin{conjecture}[AS conjecture $2$] For any $\mathbf{x} \in C_n(H^3)$,
\[ |D(\mathbf{x})| \geq 1 .\] \label{AS2} \end{conjecture}

This article is really concerned with \emph{coplanar} hyperbolic configurations. There is no loss of generality thus in assuming that
\[ \mathbf{x} \in C_n(H^2) ,\]
since one can always find a hyperbolic isometry taking a given hyperbolic plane, to another given one. However, departing from the 
standard convention of regarding $H^2$ as the subset of $H^3$ given by $z=0$, we will regard (though this is not a necessity) $H^2$ as the subset of 
$H^3$ given by $y=0$. This will make the $t_{ab}$ real, or possibly infinity (instead of making them have unit modulus, as in the case $z=0$).

Our first result can now be formulated.
\begin{theorem}[Thm. $1$] For any $\mathbf{x} \in C_4(H^2)$ such that no point in the configuration $\mathbf{x}$ lies in the convex hull of the other three, AS 
conjecture $2$ holds, namely
\[ |D(\mathbf{x})| \geq 1. \] \label{Thm1} \end{theorem}

We also reprove Y. Zhang and J. Ma's result (see \cite{MZ}), namely
\begin{theorem}[Zhang-Ma] For any $\mathbf{x} \in C_4(H^2)$ such that no three of the points in the configuration $\mathbf{x}$ are collinear, and one of the points 
in $\mathbf{x}$ lies in the convex hull of the other three, AS conjecture $1$ holds, namely
\[ D(\mathbf{x}) \neq 0 .\] \end{theorem}

Instead of defining the polynomials $p_a$ via \eqref{observer-based} (which I refer to as the ``observer''-based polynomials, and these are actually 
the Atiyah-Sutcliffe polynomials), we may alternatively define a similar, though different set of polynomials
\begin{equation} q_a(t) = \prod_{b\neq a} (t - t_{ba}) . \label{star-based} \end{equation}
We refer to these polynomials as the ``star''-based variant of the Atiyah-Sutcliffe problem. We could also define
\[ D^s(\mathbf{x}) = \frac{\det(q_1,\ldots,q_n) }{\prod_{a<b} \det(q_{ab},q_{ba})}, \]
where $q_{ab}(t) = t - t_{ba}$, up to a non-zero complex scalar factor.

One may also conjecture the analogues of the AS conjectures $1$ and $2$ for $D^s$ instead of $D$. Our final result can be formulated.
\begin{theorem}[Thm. $2$] If $\mathbf{x} \in C_n(H^2)$ is a convex coplanar configuration (by convex we mean that none of the points in the 
configuration lies in the convex hull of the other points), then
\[ D^s(\mathbf{x}) \geq \left( \prod_{k=1}^{n-1} k! \right)^{-2}. \]
This is of course enough to imply the star-based analogue of AS conjecture $1$, but not enough to imply the analogue of AS conjecture $2$, 
for such configurations. \end{theorem}

In section \ref{n4}, we prove the $n=4$ theorems: Thm. 1 and the Zhang-Ma result. Then in section \ref{general}, we prove Thm. 2.

\section{Proofs of the $n=4$ results} \label{n4}

In this section, we prove both Thm. $1$ and the Zhang-Ma theorem. We make use of computer algebra. More specifically, we make use of 
the symbolic library SymPy of Python $3$.

\begin{proof}[Proof of Thm. $1$] Let $\mathbf{x} = (\mathbf{x}_1,\ldots,\mathbf{x}_4) \in C_4(H^2)$ (our notation is explained in section \ref{intro}), such that 
none of the points in $\mathbf{x}$ lies in the convex hull of the other three. We call a coplanar configuration such that none of the points lies in the convex hull of the other 
points convex. Denote the ``circle at infinity'' $\partial H^2$ by $S^1_{\infty}$. WLOG, up to relabeling the $4$ points of the configuration (and possibly using an isometry of $H^3$ 
that preserves $H^2$ as a set, yet reverses the direction of $S^1_{\infty}$), one may assume the directions $t_{ab}$ lie on $S^1_{\infty}$ in the following 
order:
\begin{equation} t_{21},t_{31},t_{41},t_{32},t_{42},t_{12},t_{43},t_{13},t_{23},t_{14},t_{24},t_{34} . \label{roots-order} \end{equation}

The Atiyah-Sutcliffe polynomials are, in this case:
\begin{align*} p_1(t) &= (t-t_{12}) (t-t_{13}) (t-t_{14}) \\
p_2(t) &= (t-t_{21}) (t-t_{23}) (t-t_{24}) \\
p_3(t) &= (t-t_{31}) (t-t_{32}) (t-t_{34}) \\
p_4(t) &= (t-t_{41}) (t-t_{42}) (t-t_{43}) \end{align*}
The numerator $AS_4(\mathbf{x})$ of $D_4(\mathbf{x})$ is then easily computed in SymPy, in terms of the ``symbols'' $t_{ab}$.
We then define
\begin{align*} y_0 &= \sum_{a \neq b} t_{ab} \\
y_1 &= t_{31} - t_{21} \\
y_2 &= t_{41} - t_{31} \\
&\,\,\, \vdots\\
y_{11} &= t_{34} - t_{24} \end{align*}
Basically, looking at equation \eqref{roots-order}, $y_1$ is the difference between the second direction and the first, and so on, so that 
$y_k$ is the difference between the $k+1$-th direction and the $k$-th direction in the same equation, for $1 \leq k \leq 11$. We can then easily invert the above linear equations, 
and express the $t_{ab}$ in terms of the $y_k$ ($0 \leq k \leq 11$). By using the corresponding substitutions, we express $AS_4(\mathbf{x})$ in terms of the $y_k$ ($0 \leq k \leq 11$), 
instead of the $t_{ab}$. It turns out that $AS_4(\mathbf{x})$ does not depend on $y_0$ and that, moreover, it is a linear combination of monomials in the $y_k$ with positive 
coefficients only.

We then consider the denominator of $D(\mathbf{x})$, which we denote by $B_4(\mathbf{x})$. We then express the latter in terms of the $y_k$, just like we did for $AS_4(\mathbf{x})$. 
It turns out also that $B_4(\mathbf{x})$ does not depend on $y_0$ and that, moreover, it is also a linear combination of monomials in the $y_k$ with positive coefficients.

We now form the difference
$\delta_4(\mathbf{x}) = AS_4(\mathbf{x}) - B_4(\mathbf{x})$.
One may then check (using Python) that $\delta_4(\mathbf{x})$ is also a linear combination of the $y_k$ ($1 \leq k \leq 11$) with positive coefficients only. This actually proves AS conjecture $2$ in this case, namely 
for the case of $4$ distinct coplanar hyperbolic points, such that none of the points is in the convex hull of the other three. We actually obtain AS conjecture $2$ for the closure of this 
subset of $C_4(H^2)$, by a continuity argument. So we also get AS conjecture $2$ for the case where at least three of the points are collinear.
\end{proof}

\begin{proof}[Proof of the Zhang-Ma Thm.]

Our proof of the Zhang-Ma theorem is analogous to our proof of Thm. $1$. WLOG, one may assume that $\mathbf{x}_2, \mathbf{x}_3, \mathbf{x}_4$ are 
non-collinear, and that $\mathbf{x}_1$ is in the interior of the convex hull of $\mathbf{x}_2, \mathbf{x}_3, \mathbf{x}_4$. Therefore, WLOG one may assume that 
the $t_{ab}$ lie on $S^1_{\infty}$ in the following order:
\begin{equation} t_{13},t_{23},t_{21},t_{24},t_{14},t_{34},t_{31},t_{32},t_{12},t_{42},t_{41},t_{43} . \label{roots-order2} \end{equation}
We then define
\begin{align*} y_0 &= \sum_{a \neq b} t_{ab} \\
y_1 &= t_{23} - t_{13} \\
y_2 &= t_{21} - t_{23} \\
&\,\,\, \vdots\\
y_{11} &= t_{43} - t_{41} \end{align*}
Thus $y_k$ is the difference between the $k+1$-th direction and the $k$-th direction in equation \eqref{roots-order2}, for $1 \leq k \leq 11$. 
Using substitutions (arising from inverting the above linear system), one may express the numerator $AS_4(\mathbf{x})$ of $D_4(\mathbf{x})$ in terms 
of the $y_k$. It turns out that $AS_4(\mathbf{x})$ does not depend on $y_0$, and is a linear combination of monomials in the $y_k$ with only negative coefficients. 
Similarly, the denominator $B_4(\mathbf{x})$ of $D_4(\mathbf{x})$ does not depend on $y_0$ and is also a linear combination of monomials in the $y_k$ with 
only negative coefficients. Since all the $y_k$ are positive, for $1 \leq k \leq 11$, under our assumption, it follows that $D_4(\mathbf{x}) > 0$ for such configurations 
$\mathbf{x}$. This proves AS conjecture $1$ in this case, and ends our proof of the Zhang-Ma theorem.
 \end{proof}

\begin{remark} If one expresses 
\[ AS_4(\mathbf{x}) - B_4(\mathbf{x}) \]
in terms of the $y_k$, in our last proof, it turns out that it is a linear combination of the $y_k$ ($1 \leq k \leq 11$), but some coefficients are positive and some are negative. 
So, while this strategy worked fine in our proof of Thm. $1$, it unfortunately fails in the Zhang-Ma setting.
\end{remark}

\section{Proof of Thm. $2$} \label{general}

We derive an expansion of $D^s_n(\mathbf{x})$ for the star-based variant of the Atiyah-Sutcliffe determinant (defined in section \ref{intro}), which holds for any 
$\mathbf{x} \in C_n(H^3)$, from which Thm. $2$ easily follows.

The key to deriving such an expansion formula for $D^s_n$, or rather for its numerator $AS^s_n(\mathbf{x})$, is the group of symmetries of such a determinant. More 
specifically, arrange the directions $t_{ab}$ in a $2$-dimensional $n$-by-$n$ array, with nothing (say $*$) on the diagonal. Denote this array by $\mathbf{t}$. Thus
\begin{equation} \mathbf{t} = \left( \begin{array}{cccc} * & t_{12} & \ldots & t_{1n} \\
t_{21} & * & \ldots & t_{2n} \\
\vdots & \vdots & \ddots & \vdots \\
t_{n1} & t_{n2} & \ldots & * \end{array} \right) 
\end{equation} 
Let $P$ be the group of permutations of the non-diagonal entries in $\mathbf{t}$ which preserve each column of $\mathbf{t}$. Thus $P$ is a finite group of order
\[ |P| = \{(n-1)!\}^n \]
The denominator $B^s_n(\mathbf{x})$ of $AS^s_n(\mathbf{x})$ is
\begin{equation} B^s_n(\mathbf{x}) = \prod_{1 \leq a<b \leq n} (t_{ab} - t_{ba}) \end{equation}
On the other hand, the expression
\[ \sum_{\sigma \in P} (\sigma.B^s_n)(\mathbf{x}) = \sum_{\sigma \in P} \prod_{1 \leq a<b \leq n} (\sigma.t_{ab} - \sigma.t_{ba}) \]
has all the symmetry properties of $AS^s_n(\mathbf{x})$. It is polynomial in the $t_{ab}$ of degree $1$ in each variable and total degree $n$-choose-$2$, it 
is preserved by $P$ and is relatively invariant by the symmetric group $\Sigma_n$ acting on both indices of $\mathbf{t}$ simultaneously with weight the sign 
homomorphism from $\Sigma_n$ onto $\{ \pm 1 \}$ (the sign of a permutation is $1$ if it is even, and $-1$ if it is odd). A more accurate way to express the symmetry 
properties of $AS^s_n(\mathbf{x})$ is to use homogeneous coordinates to describe directions. It is then homogeneous of degree $1$ in each direction, of total degree 
$n(n-1)$, is invariant under $P$ and relatively invariant with respect to $\Sigma_n$ with respect to the sign homomorphism. The vector space of polynomials having 
these symmetries is complex $1$-dimensional (this is similar to the statement that the determinant is the unique, up to scaling, skew-symmetric multilinear function of $n$ 
vectors in $\mathbb{C}^n$). 
Hence, we can deduce that
\begin{equation} AS^s_n(\mathbf{x}) = c_n \sum_{\sigma \in P} (\sigma.B^s_n)(\mathbf{x}) \label{expansion} \end{equation}
where $c_n$ is some nonzero constant depending on $n$ only. One may then easily deduce what $c_n$ is, by evaluating for instance at collinear configurations, at which 
it is known that $D^s_n$ takes the value $1$. In fact,
\begin{equation} c_n = \frac{1}{\prod_{k=1}^{n-1} (k!)^2} .\end{equation}
This is the expansion of $AS^s_n(\mathbf{x})$ that we were seeking.

We will now apply this formula to prove Thm. $2$. Since all $n$ points in $H^2$ are such that none of the points is in the convex hull of the other points, 
one may assume WLOG that the $t_{ab}$ lie on the circle at infinity $S^1_{\infty}$ in the following order:
\begin{equation} t_{21},t_{31},\ldots, t_{n1},t_{32},t_{42},\ldots, t_{n2}, t_{12},\ldots,t_{1n},t_{2n},\ldots,t_{n-1\,n} . \label{roots-order-n} \end{equation}
But then, using formula \eqref{expansion}, we see that
\[ \delta^s(\mathbf{x}) := AS^s_n(\mathbf{x}) - c_n B^s_n(\mathbf{x}) \]
is a linear combination of nonnegative polynomials in the $t_{ab}$ with positive coefficients (all the coefficients are equal to $c_n$). This is enough to show that for such 
configurations $\mathbf{x}$, we have
\[ D^s_n(\mathbf{x}) \geq c_n ,\]
and Thm. $2$ is proved.

\acknowledge{I thank Niky Kamran, Kamal Khuri-Makdisi and Dennis Sullivan for their patience in listening to my ideas, despite not being working on 
this problem. I thank Dennis for his advice on mathematical life, in general, and I thank Niky and Kamal for encouraging me to publish my results.}

\pdfbookmark[1]{References}{ref}

\end{document}